\documentclass[12pt]{article}
\input{epsf.sty}
\usepackage{hyperref}
\usepackage{amsmath, amssymb}

\usepackage[arrow, matrix, curve]{xy}

\newtheorem {thm}{Theorem}[section]
\newtheorem {prop}[thm]{Proposition} 

\newtheorem {lem}[thm]{Lemma}

\newtheorem {defn}[thm]{Definition}

\newtheorem {cond}[thm]{Condition}

\newcommand{\qed}{\nobreak \ifvmode \relax \else
      \ifdim\lastskip<1.5em \hskip-\lastskip
      \hskip1.5em plus0em minus0.5em \fi \nobreak
      \vrule height0.75em width0.5em depth0.25em\fi}

\def\Cox{\hfill \Box}

\def\N{{\Bbb N}}

\def\Z{{\Bbb Z}}
\def\R{{\Bbb R}}

\def\e{{\varepsilon}}

\def\ba{{\backslash}}
\def\sb{{\subset}}

\def\D{\Delta}
\def\a{\alpha}

\def\ba{\setminus}

\def\d{\delta}

\def\e{\varepsilon}

\def\phi{\varphi}

\def\g{\gamma}

\def\r{\rho}

\def\s{\sigma}

\def\x{\xi}

\def\o{\omega}

\def\D{\Delta}

\def\L{\Lambda}

\def\G{\Gamma}

\def\O{{\Omega}}

\def\SS{{\Sigma}}

\def\T{\T}

\def\SS{{\cal S}}

\def\GG{{\cal G}}

\def\V|{{\Vert}}

\begin{document}
\title{Attractor properties of non-reversible dynamics w.r.t invariant Gibbs measures on the lattice}

\author{
Benedikt Jahnel
\footnote{ Ruhr-Universit\"at   Bochum, Fakult\"at f\"ur Mathematik, D44801 Bochum, Germany,
\newline
 \texttt{Benedikt.Jahnel@ruhr-uni-bochum.de}, 
\newline
\texttt{http://http://www.ruhr-uni-bochum.de/ffm/Lehrstuehle/Kuelske/jahnel.html }}
 \, and  Christof K\"ulske
\footnote{ Ruhr-Universit\"at   Bochum, Fakult\"at f\"ur Mathematik, D44801 Bochum, Germany,
\newline
\texttt{Christof.Kuelske@ruhr-uni-bochum.de}, 
\newline
\texttt{http://www.ruhr-uni-bochum.de/ffm/Lehrstuehle/Kuelske/kuelske.html
/$\sim$kuelske/ }}\, 
\,  
}

\maketitle

\begin{abstract} 
{\bf }
We consider stochastic dynamics of lattice systems with finite local state space, possibly 
at low temperature, and possibly non-reversible. 
We assume the additional regularity properties on the dynamics: 

a) There is at least one stationary measure which is a Gibbs measure for an absolutely summable potential $\Phi$.  

b) Zero loss of relative entropy density under dynamics implies the Gibbs property with the same $\Phi$. 

We prove results on the attractor property of the set of Gibbs measures for $\Phi$: 
%
%
%
%

1. The set of weak limit points of any trajectory of translation-invariant measures contains at least one Gibbs state for $\Phi$. 

2. We show that if all elements of a weakly convergent sequence of measures are Gibbs measures for a sequence of some translation-invariant summable potentials with uniform bound, then the limiting measure must be a Gibbs measure for $\Phi$.

3. We give an extension of the second result  to trajectories which are allowed to be non-Gibbs, 
but have a property of asymptotic smallness of discontinuities.  An example for this situation is the time evolution from a low temperature Ising measure by weakly dependent spin flips. 
\end{abstract}

\smallskip
\noindent {\bf AMS 2000 subject classification:} 82C20,
82C05, 60K35.

 \smallskip
\noindent {\bf Keywords:} Markov chain, PCA, IPS, non-equilibrium, non-reversibility, attractor property, relative entropy, Gibbsianness, non-Gibbsianness, synchronisation.

\vfill\eject

\vfill\eject

\newpage 
\section{Introduction} 
The study of non-equilibrium statistical mechanics models away from their time-stationary equilibrium states and their relaxation or non-relaxation into equilibrium is an active field of research in the theoretical physics community as well as in the mathematics of Markov processes \cite{Li85,Ma06}. If there exists initial data that does not converge into the equilibrium state (even in the presence of a unique time-stationary measure) the model is called non-ergodic and examples can be found in  \cite{JaKu12,JaKu14b,MaSh11,ChMa11}. If there is relaxation of all initial data towards some set of measures we call this set an attractor. 
In certain settings when there exists a periodic orbit of measures, this phenomenon is also called synchronisation and represents a common feature in many areas of science and engineering. Examples are found experimentally and in simulations in
the study of neuronal pulses of the brain or digital communications receivers, and partially understood theoretically, mostly in mean field like the Kuramoto model (see e.g. references \cite{AcBoPeRiSp05,GPP12,DH96,GiPaPePo12,PiRoKu01}.)

\medskip 
The purpose of this note is to provide some criteria which allow to control the approach to attractors beyond situations with 
weak interactions and beyond reversible dynamics. The criteria will be formulated in terms of regularity of trajectories, in a sense to 
be described below. 

\medskip
Restricting to translation-invariant statistical mechanics models on the lattice makes available the powerful relative entropy techniques \cite{Ma06,EnFeSo93} highlighted for example already in the Gibbs variational principle \cite{Ge11}. The main idea for the dynamical models is to look at the change of the relative entropy density of a given measure w.r.t a time-stationary measure under the evolution. It turns out, that this change is non-positive under rather general assumptions \cite{EnFeSo93}.
The use of the relative entropy density as a Lyapunov function is subtle because it is not a weakly continuous functional in the space of measures. 
More work is needed, requiring some regularity of the time-evolved measures. 
 Notice also that the relative entropy density can not distinguish between different Gibbs measures for the same potential. So, 
 in the presence of phase-transition of the equilibrium model, that is 
 when there are more than one Gibbs measures corresponding 
to the potential of the time-stationary measure, the entropy method can at most ensure attraction of the whole 
set of Gibbs measures.  

\medskip
The relative entropy approach has been used to prove that measures having zero entropy loss under the dynamics w.r.t a time-stationary Gibbs measure are Gibbs measures for the same potential in the examples of \cite{Ho71,HiSh75,DaLoRo02,HoSt77}. These concern stochastic Ising models but also more general 
probabilistic cellular automata without reversibility assumption where the aspect of attractivity from initial states away from the invariant set was not discussed. 
In this note we provide results on the limiting behavior of trajectories with general initial data for general translation-invariant discrete-time Markov processes (DMP) and continuous-time interacting particle systems (IPS)  on $\{1,\dots,q\}^{\Z^d}$ assumed to have the above zero entropy loss property. The previous examples show that this hypothesis is satisfied in a number 
of important cases. Let us also mention the case of the well known symmetric exclusion process (SEP)(see for example \cite{Li85} Chapter VIII). Here the stationary measures $\mu$ are product measures and zero entropy loss w.r.t to $\mu$ implies that the time-evolved measure is a mixture of product measure.

\medskip
One more specific motivation for this note comes from the investigation of a class of non-equilibrium statistical mechanics models with $d\geq 3$ proposed by the authors in \cite{JaKu12,JaKu14b}. Here the dynamics is given by a non-reversible probabilistic cellular automaton (PCA) with exponentially localized updating rule (see \cite{JaKu14b}) respectively by an IPS dynamics (see \cite{JaKu12}). Both have the property to create non-trivial periodic orbits of extremal translation-invariant Gibbs measures. The primary focus in \cite{JaKu12,JaKu14b} was to demonstrate that there can be models showing non-ergodic behavior in the presence of a unique translation-invariant time-stationary measure. Once this is established, it is natural to ask if and how the dynamics drives any starting measure into the periodic orbit. That this is indeed the case for a mean-field version of the IPS dynamics is one of the main results in \cite{JaKu14}. Let us mention again the Kuramoto model which also is a mean-field statistical mechanics system driven by its Langevin dynamics. Here similar results have been obtained see \cite{AcBoPeRiSp05,GPP12,DH96,GiPaPePo12}.  

\subsection{Strategy and main results}
The main objective of the present paper is to give criteria for 
a given set of measures, containing at least one invariant Gibbs measure w.r.t PCA and IPS dynamics, to be an attractor for a stochastic dynamics in a lattice setup. 
Let us mention that if the dynamics has specific monotonicity properties like "attractivity" (in the sense 
of stochastic domination being preserved by dynamics) coupling arguments can be used to derive attractor properties, see \cite{Li85} Chapter III Section 2. 
Here we want to treat cases also beyond that. 

\medskip
The strategy is exemplified in the very special case of the stochastic Ising model (also called Glauber dynamics) for a not necessarily ferromagnetic translation-invariant Hamiltonian with local spin space $\{\pm 1\}$
and finite range interactions by Holley \cite{Ho71} and for not necessarily finite range but fast decaying interactions by Higuchi and Shiga \cite{HiSh75}.
Here it has been proved that any limit measure 
of a sequences of measures (propagated by the Glauber dynamics) must be a Gibbs measure. 
The main tool in both cases is to consider translation-invariant measures and the change in relative entropy density between those and the Gibbs measures under the dynamics.
To be more precise, the strategy is as follows. First it is shown, that the time-derivative of the relative entropy density between measures away from the Gibbs measure and the Gibbs measures is non-increasing under time-evolution, i.e $g(\nu_t|\mu):=\frac{d}{ds}_{|s=t}h(\nu_s|\mu)\leq0$ and thus since $\infty>h(\nu_t|\mu)\geq0$, $\lim_{t\uparrow\infty}g(\nu_t|\mu)=0$. This fact although is not the crucial point since it is true for rather general transformations of measures as mentioned above. What is important is to prove that while the relative entropy density itself is semicontinuous also the time-derivative of the relative entropy density is semicontinuous in the useful direction,
more precisely one has upper semicontinuity of $\nu_t\mapsto g(\nu_t|\mu)$. 
The semicontinuity from above guarantees that for a convergent sequence of measures where the sequence of time-derivatives of the relative entropy densities goes to zero, also for the limiting measure the time-derivative of the relative entropy density is zero, i.e. for $\lim_{t\uparrow\infty}\nu_t=\nu_*$ in the weak sense we have
$$0=\lim_{t\uparrow\infty}g(\nu_t|\mu)=g(\nu_*|\mu).$$
The final step of the proof is often referred to as "Holley's argument" which uses the zero entropy loss property of the limiting measures $g_L(\nu_*|\mu)=0$ to show a single-site DLR equality for $\nu_*$. In other words 
for the stochastic Ising model any measure where the time-derivative of the relative entropy density is zero has to be a Gibbs measure.

\medskip

Investigating conditions under which the possible discontinuity of $g(\nu_t|\mu)$ for general models can be beaten, we arrive at the following results, 
assuming a Holley regularity condition (Condition \ref{zero entropy loss condition}).

\medskip
In our first main result,  Theorem \ref{ExistenceGibbsLimits}, we show that in both cases, discrete-time and continuous-time stochastic dynamics, at least one weak $\o$-limit point (a cluster point of the trajectory in infinite time) has to be a Gibbs measure for the same potential as the 
time-stationary measure. 

In our second main result, Theorem \ref{CorAttractivity}, 
we show that if all elements of a weakly convergent sequence are Gibbs measures for a uniformly bounded sequence of some translation-invariant summable potentials, 
which means that no Gibbsian pathologies persist along the trajectory for large times, then the limiting measure must be a Gibbs measure for the same potential as the given time-stationary Gibbs measure. 

In our final result, Theorem \ref{MainTheorem},
we show, that in case of the continuous-time dynamics the second result holds under weaker conditions. 
The Gibbsianness assumption on the trajectory may be replaced by a uniform non-nullness condition together with martingale convergence of single-site conditional probabilities uniformly in the trajectory. This can be seen as a property of asymptotic smallness of non-Gibbsian pathologies under time-evolution.
The proof is based on a representation of the relative entropy loss (valid for non-null probability measures) we derive in Proposition \ref{EntropyLossGen}.
An explicit example for a sequence of time-evolved measures which are non-Gibbs for all sufficiently large times but satisfy the conditions 
is the initial low-temperature Ising model in zero field under infinite-temperature Ising dynamics (see \cite{EnFeHoRe02}).


\bigskip

\section{Entropy decay under time-discrete and time-continuous interacting systems}

Consider translation-invariant probability measures $\mu$ and $\nu$ on the configuration space $\{1,\dots,q\}^{\Z^d}$ equipped with the usual product topology and the Borel sigma-algebra. For a finite set of sites $\L\subset\Z^d$ define the \textit{local relative entropy} via
\begin{equation*}\label{LRE}
\begin{split}
h_\L(\nu|\mu):=\sum_{\o_\L\in\{1,\dots,q\}^\L}\nu({\o_\L})\log\frac{\nu({\o_\L})}{\mu({\o_\L})}.
\end{split}
\end{equation*}
and the \textit{relative entropy density} via
\begin{equation*}\label{SRE}
\begin{split}
h(\nu|\mu):=\lim_{\L\uparrow\Z^d}\frac{1}{|\L|}h_\L(\nu|\mu)
\end{split}
\end{equation*}
where $\L$ runs over hypercubes centered at the origin, whenever the limit exists. We use notations like $\o_\L:=\{\s\in\{1,\dots,q\}^{\Z^d}: 1_{\o_\L}(\s)=1\}$, $\o_\D\o_{\L\setminus\D}:=\o_\D\cap\o_{\L\setminus\D}$, $\D^c:=\Z^d\setminus\D$ etc.

\bigskip
Further consider two types of translation-invariant Markovian dynamics on 
$\{1,\dots,q\}^{\Z^d}$:
\begin{enumerate}
\item Discrete-time Markov Processes (DMP) characterized by time-homogeneous transition kernels $P(\s,\cdot)$ which are also assumed to be continuous in the first entry w.r.t the product topology. 
\item Interacting particle systems (IPS) characterized by time-homogeneous generators $L$ with domain $D(L)$ and its associated Markovian semigroup $(P^L_t)_{t\geq0}$.
\end{enumerate}
Standard examples of DMP are the so-called (strict) \textit{probabilistic cellular automata} (PCA) characterized by the fact, that the transition kernels factorize, i.e. $P(\s,\eta_\L)=\prod_{i\in\L}P_i(\s,\eta_i)$, see \cite{DaLoRo02}. Also more general \textit{PCA with exponentially localized update kernel} can be considered, see for example \cite{JaKu14b,Ku84}. For the IPS we adopt the exposition given in \cite{Li85} Chapter I: In all generality we let the generator $L$ be given via jump-measures $c_\D(\eta,d\xi_\D)$ in finite volumes $\D\subset\Z^d$, continuous in the starting configurations $\eta\in\{1,\dots,q\}^{\Z^d}$
\begin{equation*}
Lf(\eta)=\sum_\D\int_{\{\xi:\xi_{\D^c}=\eta_{\D^c}\}}c_\D(\eta,d\xi)[f(\xi)-f(\eta)]
\end{equation*}
where the summation is over all finite sets of sites and $f\in D(L)$. To ensure well-definedness, the jump-measures must satisfy a number of conditions, most importantly the single-site jump-intensities have to be bounded, i.e for
$c_\D:=\sup_{\eta}c_\D(\eta,\{1,\dots,q\}^\D)$
we assume $\sum_{\D\ni0}c_\D<\infty$.

\bigskip
The relative entropy density can be understood as a measure of closeness between the probability measure in the first and second entry. Accordingly the change in relative entropy density under the application of the dynamics measures the change in distance between the two probability measures. Let us recall some important facts about Gibbs measures and relative entropy densities.

\begin{lem}\label{EnterLemma}
Let $(\O^{\Z^d},\SS)$ and $(\tilde\O^{\Z^d},\tilde\SS)$ be measurable spaces of lattice configurations and $T$ any translation-invariant probability kernel from $(\O^{\Z^d},\SS)$ to $(\tilde\O^{\Z^d},\tilde\SS)$, i.e for all $i\in\Z^d$, $\tilde A\in \tilde\SS$ and $\eta\in\O^{\Z^d}$ we have $T(\tilde A|\eta)=T(\tilde A_{\theta(i)}|\eta_{\theta(i)})$ where $\tilde A_{\theta(i)}$ denotes the lattice translates of $\tilde A$ by $i$ (respectively $\eta_{\theta(i)}$ the translate of $\eta$ by $i$). Then $h(T\nu|T\mu)\leq h(\nu|\mu)$ for all translation-invariant probability measures $\nu,\mu$ on  $(\O^{\Z^d},\SS)$.
\end{lem}
For the proof see for example \cite{EnFeSo93} Lemma 3.3.
\begin{lem}\label{RED_Exists}
Let $\nu$ and $\mu$ be translation-invariant measures on the measurable configuration space $(\O,\SS)$ and $\mu$ a Gibbs measure for the Gibbsian specification $\g^\Phi$. Then the relative entropy density $h(\nu|\mu)$ exists and depends only on $\nu$ and $\Phi$.
\end{lem}
For details on Gibbs measures and their definition via the DLR equation for models given in terms of Gibbsian specifications $\g^\Phi$ see \cite{Ge11} Chapter 1 and 2. The lemma is part of Theorem 15.30 in \cite{Ge11}. The \textit{Gibbs variational principle} states that under the conditions of the preceding lemma $h(\nu|\mu)=0$ if and only if $\nu$ is a Gibbs measure for the Gibbsian specification $\g^\Phi$. Note, that for the existence of $h$, the requirement of $\mu$ to be a Gibbs measure can be relaxed considerably. The appropriate notion is that of \textit{asymptotically decoupled measures} as defined in \cite{Pf02,KuLeRe04}.

\bigskip
Consider a 
model given in terms of the Gibbsian specification $\g^\Phi$ and a translation-invariant DMP or IPS dynamics. 
Let us assume that for the dynamics the following \textit{zero entropy loss condition} holds:
\begin{cond}\label{zero entropy loss condition}
There exists a translation-invariant and time-stationary Gibbs measure $\mu$ for $\g^\Phi$. Further, for any translation-invariant measure $\nu$ with
\begin{enumerate}
\item $g_P(\nu|\mu):=h(P\nu|\mu)-h(\nu|\mu)=0$ it follows that $\nu$ is a Gibbs measure for $\g^\Phi$ (in the case of discrete-time dynamics),
\item $g_L(\nu|\mu):=\lim_{\L\uparrow\Z^d}\frac{1}{|\L|}\frac{d}{dt}_{|t=0}h_\L(P_t^L\nu|\mu)=0$ it follows that $\nu$ is a Gibbs measure for $\g^\Phi$ (in the case of continuous-time dynamics).
\end{enumerate}
\end{cond}
Such a condition is proved to hold in continuous time for example for the stochastic Ising model \cite{Li85,Ho71,HoSt77} or more general Glauber dynamics and even non-reversible dynamics see \cite{JaKu12}. In discrete time examples are given in \cite{DaLoRo02,JaKu14b}.

\bigskip
\textbf{Remark: }  We provide another example where zero entropy loss implies Gibbsianness w.r.t. to the same potential 
as the reference measure in the second slot, however after also taking into account a global preservation of density of particles 
which is conserved by the dynamics. 

Let us consider the above condition for the well known symmetric exclusion process (SEP) on the $d$-dimensional integer lattice (see for example \cite{Li85} Chapter VIII)
\begin{equation*}
Lf(\eta):=\sum_{x\in\Z^d}\sum_{y:y\sim x}\eta(x)(1-\eta(y))[f(\eta^{xy})-f(\eta)]
\end{equation*}
where $y\sim x$ denotes nearest neighbors relation of $x$ and $y$, $\eta^{x,y}$ stands for the configuration equal to $\eta$ except for the sites $x$ and $y$ where 
it is flipped.
$f$ is a sufficiently smooth observable.

It is known that the extremal stationary measures are the product measures $\mu_{\r}$ and a classification of their 
basins of attraction in terms of densities of the initial measure can be given (\cite{Li85} Chapter VIII Theorem 1.47). 
A translation-invariant ergodic (that is tail-trivial) initial measure with density $\r$ converges to $\mu_{\r}$. From this it is clear that 
the limit of any translation-invariant initial measure is the corresponding mixture over product measures. 
Product measures are Gibbs measures without interaction, and product measures with different densities 
are Gibbs measures for different specifications.
On the one hand product measures are simpler than the Gibbs measures with interaction and their possible phase transitions
we have encountered in our other examples. 
On the other hand the SEP is more general than our other examples 
since possible limits correspond to sets of specifications and not a single specification. 
Let us see that our condition is consistent with this picture
by showing that there are no other ergodic measures with fixed density 
which have zero entropic loss w.r.t. to one of the  $\mu_{\r}$'s.

%

Equating the entropy loss of a translation-invariant measure $\nu$ w.r.t one of the invariant product-measures $\mu_{\r}$ in this case to zero, 
we immediately see that the dependence on $\r\in[0,1]$ drops out. 
Indeed, by translation invariance
\begin{equation*}\label{SEP}
\begin{split}
g_{SEP}&(\nu|\mu_\r)=\sum_{i\sim 0}\int\nu(d\eta)\nu(\s_0=1\,\,\s_i=0|\eta_{\{0,i\}^c})\log\frac{\nu(\s_0=0\,\,\s_i=1|\eta_{\{0,i\}^c})}{\nu(\s_0=1\,\,\s_i=0|\eta_{\{0,i\}^c})}\cr
&=\sum_{i=1}^d\int\nu(d\eta)[\nu(1_0\,\,0_{e_i}|\eta_{\{0,e_i\}^c})-\nu(0_0\,\,1_{e_i}|\eta_{\{0,e_i\}^c})]\log\frac{\nu(0_0\,\,1_{e_i}|\eta_{\{0,e_i\}^c})}{\nu(1_0\,\,0_{e_i}|\eta_{\{0,e_i\}^c})}.\cr
\end{split}
\end{equation*}
This implies $g_{SEP}(\nu|\mu_\r)\leq0$ and if we set $g_{SEP}(\nu|\mu_\r)=0$ we have
\begin{equation}\label{SEP2}
\begin{split}
\nu(1_0\,\,0_{e_i}|\eta_{\{0,e_i\}^c})=\nu(0_0\,\,1_{e_i}|\eta_{\{0,e_i\}^c})
\end{split}
\end{equation}
for $\nu-a.a.$ $\eta$ and $i\in\{1,\dots,d\}$. This implies $\nu(\eta_V|\eta_{V^c})=\nu(\pi_V(\eta_V)|\eta_{V^c})$ for any finite volume $V$ and any permutation $\pi_V(\eta_V)$ of the finite configuration $\eta_V$. Indeed, we can assume $V$ to be a box since there exists a box $B\supset V$ and if we assume $\nu(\eta_B|\eta_{B^c})=\nu(\pi_B(\eta_B)|\eta_{B^c})$ for any permutation $\pi_B$, of course also $\nu(\eta_B|\eta_{B^c})=\nu(\pi_V(\eta_B)|\eta_{B^c})$ and thus
\begin{equation*}
\begin{split}
\nu(\pi_V(\eta_V)|\eta_{V^c})=\frac{\nu(\pi_V(\eta_B)|\eta_{B^c})}{\nu(\eta_{B\setminus V}|\eta_{B^c})}=\frac{\nu(\eta_B|\eta_{B^c})}{\nu(\eta_{B\setminus V}|\eta_{B^c})}=\nu(\eta_V|\eta_{V^c}).
\end{split}
\end{equation*}
Further, any finite permutation $\s$ can be realized as a finite product of nearest-neighbor transpositions $\pi_{i,j}(\eta_B)=\eta_{B\setminus\{i,j\}}(\eta_{j})_i(\eta_{i})_{j}$ where $j\sim i$. If $\eta_i=\eta_j$, there is nothing to show. If $\eta_i\neq\eta_j$ by \eqref{SEP2}, translation-invariance and the elementary definition of conditional probability we have
\begin{equation*}
\begin{split}
\nu(\pi_{i,j}(\eta_B)|\eta_{B^c})&=\nu(\pi_{i,j}(\eta_{\{i,j\}})|\eta_{\{i,j\}^c})\nu(\eta_{B\setminus\{i,j\}})|\eta_{B^c})\cr
&=\nu(\eta_{\{i,j\}})|\eta_{\{i,j\}^c})\nu(\eta_{B\setminus\{i,j\}})|\eta_{B^c})=\nu(\eta_B|\eta_{B^c}).\cr
\end{split}
\end{equation*}
From the invariance of the conditional probabilities w.r.t finite permutations, it follows, that $\nu$ is invariant w.r.t finite permutations, in other words 
exchangeable. By de Finetti's Theorem (see \cite{Ge11} Example 7.16 and 7.31) it is thus a mixture of product measures $\nu=\int_0^1\mu_\r m_\nu(d\r)$ where $m_\nu(d\r)$ is a unique probability measure on the product-measures together with the evaluation sigma-algebra (for more details see \cite{Ge11} Chapter 7: Extreme decomposition). 
Since the only tail-trivial mixtures of product measures are the pure product measures themselves, our claim follows.


%

\bigskip
Let us state our first result about attractor properties.

\begin{thm}\label{ExistenceGibbsLimits}
Assume Condition \ref{zero entropy loss condition} holds with Gibbs measure $\mu$ for $\g^\Phi$. Let $\nu_0$ be any translation-invariant starting measure. Then the set $C$ of all weak limit points 
of the sequence $\nu_n:=P^n\nu_0$ respectively $\nu_{t}:=P_{t}^L\nu_0$  
contains translation-invariant Gibbs measures for $\g^\Phi$. 
\end{thm}

\textbf{Proof: }First we note that the set $C$ of all weak limit points is weakly compact. 
Further, the map $\nu\mapsto h(\nu|\mu)$ is lower semicontinuous by \cite{Ge11} Theorem 15.39 and 
hence the infimum of 
$\nu\mapsto h(\nu|\mu)$ as a map from $C$ to $\R^+_0\cup\{+\infty\}$ is attained in some $\nu_{*}\in C$.

Suppose $h(\nu_{*}|\mu)>0$, then $\nu_{*}$ in that case is not a Gibbs measure for the same potential as $\mu$ by the Gibbs variational principle (see \cite{Ge11} Theorem 15.39). Further by Condition \ref{zero entropy loss condition} for the discrete-time case $$h(P\nu_{*}|\mu)<h(\nu_{*}|\mu)$$ 
and for the continuous-time case for all $t>0$ $$h(P_t^L\nu_{*}|\mu)<h(\nu_{*}|\mu).$$
But this is a contradiction since $P\nu_{*}$ respectively $P_t^L\nu_{*}$ are also weak limit points by the continuity of $P$ and $P_t^L$.
$\Cox$

\bigskip
The preceding theorem in particular implies that for convergent trajectories the then unique $\o$-limit measure (the then unique cluster point of the trajectory in infinite time) must be a Gibbs measure for $\g^\Phi$. Under Condition \ref{zero entropy loss condition} this also follows from the fact that the limiting measure is invariant for the dynamics (see \cite{Li85} Proposition 1.8. for the IPS case, the DMP case follows easily by the same arguments).

\subsection{Attractor properties along Gibbsian trajectories}
The next result makes the assumption that all but finitely many elements of the converging subsequence are translation-invariant Gibbs measures for a uniformly bounded sequence of translation-invariant potentials. Here we define the norm
$\Vert\Phi\Vert:=\sum_{A\ni 0}\Vert\Phi_A\Vert_\infty$. As we will see, the benefit from this is the fact, that the change of entropy as a function of the first entry $\nu\mapsto g_P(\nu|\mu)$ and $\nu\mapsto g_L(\nu|\mu)$ is continuous along such a sequence of measures. Let us note that for the attractor property of the set of Gibbs measures, checking upper semicontinuity of the change of the relative entropy density would be sufficient, see also \eqref{UpperSemi}. This is what is in fact done in \cite{Ho71,HiSh75}.

\begin{thm}\label{CorAttractivity} 
Assume Condition \ref{zero entropy loss condition} holds with Gibbs measure $\mu$ for $\g^\Phi$. Let $\nu_0$ denote an arbitrary translation-invariant starting measure. Further let $(\nu_{n_k})_{k\in\N}$ (resp. $(\nu_{t_k})_{k\in\N}$) be any weakly convergent subsequence of the sequence of time-evolved measures $\nu_n:=P^n\nu_0$ (resp. $\nu_t:=P_t\nu_0$ with $t_k\uparrow\infty$) and let $\nu_*$ denote its weak limit. 
Suppose that 
\begin{enumerate}
\item for all $n_k$ (resp. $t_k$), the measures $\nu_{n_k}$ and $\nu_{n_k+1}$ (resp. $\nu_{t_k}$) are Gibbs measures for some translation-invariant potentials $\Phi_{n_k}$ and $\Phi_{n_k+1}$ (resp. $\Phi_{t_k}$),
\item the sequences of potentials $(\Phi_{n_k})_{k\in\N}$  and $(\Phi_{n_k+1})_{k\in\N}$ (resp. $\Phi_{t_k}$) are uniformly bounded, i.e. $S:=\sup_{k}\max\{\Vert\Phi_{n_k}\Vert,\Vert\Phi_{n_k+1}\Vert\}<\infty$ (resp. $S:=\sup_{k}\Vert\Phi_{t_k}\Vert<\infty$).
\end{enumerate}
Then, necessarily $\nu_*$ also is a Gibbs measure for $\g^\Phi$. 
\end{thm}

\textbf{Remark: }
1. Notice that the map between potentials and Gibbsian specifications is one-to-one when the equivalence relation of physical equivalence, \cite{Ge11,EnFeSo93}, is used on the space of potentials. For more details on the relation of specifications and potentials see \cite{Su73,Ko74}, in particular for the regrouping of potentials see \cite{Ku01}.
Hence one wants to exploit the theorem for {\em useful choices} of representatives in the class of physically equivalent potentials.
This is the same as looking at equivalence classes of 
physically equivalent potentials in the definition of the Banach space of potentials. In that sense we also prove, that $\lim_{k\to\infty}\Phi_{n_k}=\Phi_*$ (resp. $\lim_{k\to\infty}\Phi_{t_k}=\Phi_*$) exists with $\Vert\Phi_{*}\Vert\leq S$ and $\g^\Phi=\g^{\Phi_*}$.

2. However, having said this, the property of $\nu_{n_k}$ (resp. $\nu_{t_k}$) being Gibbs may depend strongly on the starting measure, and we can not expect it to be true universally, given the many examples of non-Gibbsian measures known to appear 
under time-evolutions \cite{KuRe06,EnRu09,EnFeHoRe02,EnFeHoRe10,ErKu10,FeHoMa13b,FeHoMa13a,HoReZu13,KuOp08,KuNy07}. This is the reason 
for our desire to relax the hypothesis and include cases of non-Gibbsian behavior, see below. 

\bigskip

{\bf Proof of Theorem \ref{CorAttractivity}: }First notice, that also the sequence $(P\nu_{n_k})_{k\in\N}$ is weakly convergent with $\lim_{k\to\infty}P\nu_{n_k}=P\nu_*$ since $P$ is continuous. 

\medskip
\textbf{Step 1: }
In order to see that $(\Phi_{n_k})_{k\in\N}$ is a convergent sequence, we show that $(\Phi_{n_k})_{k\in\N}$ is a Cauchy sequence in the Banach space of Gibbs potentials with norm\linebreak $\Vert\Phi\Vert_0:=\sum_{A\ni0}|A|^{-1}\Vert\Phi_A\Vert_\infty$ modulo physical equivalence. By \cite{EnFeSo93} formula (2.65) we can recover the corresponding potentials in the sense that 

\begin{equation*}\label{EnFeSo-2.65}
\begin{split}
\Vert\Phi_1-\Phi_2\Vert_0=\frac{1}{|\L|}\Big\Vert\log\frac{d\nu_1|_\L}{d\nu_2|_\L}\Big\Vert_C-\frac{o(|\L|)}{|\L|}
\end{split}
\end{equation*}
where also
$\Vert f\Vert_C:=\sup_{c: \text{constant}}\Vert f-c\Vert_\infty$. 
The error term may a priori depend on the potentials $\Phi_1$ and $\Phi_2$. 
By the uniform boundedness of the sequence of potentials, however the error term can be bounded by the uniform expression $\textit{const }S\frac{|\partial\L|}{|\L|}$ 
where $\partial\L$ denotes the boundary of $\L$. Let $\e>0$ and choose a centered cube $\L$ such that $\textit{const }S\frac{|\partial\L|}{|\L|}<\frac{\e}{2}$. Further by the weak convergence of the measures there exists $N_{\e}$ such that 
\begin{equation*}\label{EnFeSo-2.65weiter}
\begin{split}
\frac{1}{|\L|}\Big\Vert\log\frac{d\nu_{n_s}|_{\L}}{d\nu_{n_t}|_{\L}}\Big\Vert_C<\frac{\e}{2}
\end{split}
\end{equation*}
for all $s,t\geq N_{\e}$ using also the uniform non-nullness of all measures in the trajectory. Uniform non-nullness follows easily from the uniform boundedness of the potentials. Consequently for all $s,t\geq N_{\e}$
\begin{equation*}\label{EnFeSo-2.65weiter}
\begin{split}
\Vert\Phi_{n_s}-\Phi_{n_t}\Vert_0\leq\frac{1}{|\L|}\Big\Vert\log\frac{d\nu_{n_s}|_{\L}}{d\nu_{n_t}|_{\L}}\Big\Vert_C+\textit{const }S\frac{|\partial\L|}{|\L|}<\e.
\end{split}
\end{equation*}
Notice that for the limiting potential we also have $\Vert\Phi_*\Vert\leq S$:
Indeed, if we assume $\Vert\Phi_*\Vert\geq S+\e$ for some $\e>0$ then there exists $N\in\N$ such that $\sum_{A\ni0,|A|\leq N}\Vert\Phi_{A,n_k}-\Phi_{A,*}\Vert_\infty>\frac{\e}{2}$ for all $k\in\N$. But we have $$\lim_{k\to\infty}\sum_{A\ni0,|A|\leq N}\frac{1}{|A|}\Vert\Phi_{A,n_k}-\Phi_{A,*}\Vert_\infty=0$$ 
for all $N$, a contradiction.
Replacing $n_k$ by $t_k$ we get the same result for the continuous-time case.

\medskip
\textbf{Step 2: }For any translation-invariant starting measure $\nu_0$ we have that $h(\nu_n|\mu)$ (resp. $h(\nu_t|\mu)$) is a non-increasing sequence of non-negative numbers. (Note that the relative entropy density is smaller than infinity, due to the finite local state space and since the measure in the second slot is a Gibbs measure.) 

This sequence hence has a limit (which may a priori be strictly bigger than zero), but  from this follows that the sequence of entropy losses $g(\nu_n|\mu)$ (we write $g(\nu_n|\mu)$ for both $g_P(\nu_n|\mu)$ and $g_L(\nu_{t_n}|\mu)$)
converges to zero. 
We would like to conclude that from $\lim_{n\uparrow\infty} \nu_n= \nu_*$ in a weak sense and $\lim_{n\uparrow\infty} g(\nu_n|\mu)=0$ it  follows that $g(\nu_*|\mu)=0$. Then we know that $\nu_*$ has to be Gibbs for $\g^\Phi$ by Condition \ref{zero entropy loss condition}.


\medskip
\textbf{The discrete-time case: }
%
%
Now suppose that  $\nu$ is a Gibbs measure for some translation-invariant potential $\Phi_\nu$ and $P\nu$ is a Gibbs measure for some translation-invariant potential $\Phi_{P\nu}$ and $\mu$ is a Gibbs measure for some translation-invariant potential $\Phi$. 
We use the decomposition of the relative entropy as in 
\cite{Ge11} formula (15.32) into the pressure $p$ of the potential 
for the measure in the second slot, the expectation $\langle\cdot,\cdot\rangle$ of the local energy density of the potential of the measure in the second slot w.r.t the first measure, and the relative entropy density of the first measure, i.e.
$$h(\nu|\mu)=p(\Phi)+ \langle \nu, \Phi\rangle + h(\nu|u)$$ 
with $\langle \nu, \Phi\rangle:=\nu(\sum_{A\ni0}|A|^{-1}\Phi_A)$
and $p(\Phi):=\lim_{\L\uparrow G}|\L|^{-1}\log Z_\L^{\Phi}(\o)$ where $Z_\L^{\Phi}(\o)$ is the partition function for $\Phi$ evaluated at some arbitrary boundary condition $\o$ outside $\L$.  
We use this to write the entropy loss as
\begin{equation*}\label{Attractivity-b}
\begin{split}
g_P(\nu|\mu)&=p(\Phi)+ \langle P\nu, \Phi\rangle + h(P\nu|u)
-[p(\Phi)+ \langle \nu, \Phi\rangle + h(\nu|u)]\cr
&=\langle P\nu-\nu, \Phi\rangle+h(P\nu|u)- h(\nu|u).
\end{split}
\end{equation*}
The first term is weakly continuous in $\nu$ and causes no problems since also $P$ is continuous. 
For the second term, a priori, we have no knowledge about interchangeability of limits. Another way of considering this issue is to rewrite the relative entropy density $h(\nu|u)=\nu(\tilde h)$ as an $\nu$-expectation of a certain function $\tilde h$ as in \cite{Ge11} Theorem 15.20 where the function $\tilde h$ is not quasilocal but tail measurable. Hence convergence of expected values w.r.t a locally convergent sequence of measures is not guaranteed. 
In fact if the identity $$\lim_{k\to\infty}[h(P\nu_{n_k}|u)-h(\nu_{n_k}|u)]=h(P(\lim_{k\to\infty}\nu_{n_k})|u)-h(\lim_{k\to\infty}\nu_{n_k}|u)$$ 
were true, the result would follow. As we will show now, the uniform Gibbsianness assumption on the trajectory is sufficient to ensure such an identity.
The difference in specific entropies, assuming Gibbsianness of the two measures, can be written as
\begin{equation*}\label{Attractivity-b}
\begin{split}
h(P\nu|u)-h(\nu|u)&=\langle\nu,\Phi_\nu\rangle-\langle P\nu,\Phi_{P\nu}\rangle+p(\Phi_\nu)-p(\Phi_{P\nu}).\cr
\end{split}
\end{equation*}
The specific energy $\nu,\Phi\mapsto\langle\nu,\Phi\rangle$ is jointly continuous w.r.t the weak topology for the probability measures and the topology of convergence for the potentials (see \cite{Ge11} Remark 15.26 (2)). The same argument applies for the second term on the r.h.s of the last display. By the first part of the proof the potentials are in fact convergent and thus one can deduce interchangeability of limits. The pressure terms are continuous as functions of the potentials in the topology of uniform convergence generated by $\Vert\cdot\Vert$ (see \cite{EnFeSo93} Proposition 2.58 (b) and Proposition 2.56 (d)) and hence limits in the entropic loss can be interchanged.

\medskip
\textbf{The Continuous-time case: }We need to show for a Gibbsian sequence $(\nu_{t_k})_{k\in\N}$ that
$$\lim_{k\to\infty}g_L(\nu_{t_k}|\mu)=g_L(\lim_{k\to\infty}\nu_{t_k}|\mu).$$ 
In what follows the representation of the entropy loss $g_L$ in terms of the pairing given in \eqref{GibbsRepr} will be important. To derive this representation let us split $g_L$
into several parts. We have
\begin{equation}\label{BeforeLimit}
\begin{split}
\frac{d}{dt}_{|t=0}h_\L(P_t^L\nu|\mu)=\sum_{\o_\L}\nu(L1_{\o_\L})\log\nu(\o_\L)-\sum_{\o_\L}\nu(L1_{\o_\L})\log\mu(\o_\L).
\end{split}
\end{equation}
By properties of the relative entropy density, namely Lemma 15.28 in \cite{Ge11} and the Gibbsianness of the measures involved we can for the r.h.s of \eqref{BeforeLimit} also consider 
\begin{equation}\label{BeforeLimit2}
\begin{split}
\sum_{\o_\L}\nu(L1_{\o_\L})H_\L(\o_\L\xi_{\L^c})-\sum_{\o_\L}\nu(L1_{\o_\L})H_\L^\nu(\o_\L\xi_{\L^c})
\end{split}
\end{equation}  
and the error is of boundary order. Here $H$ and $H^\nu$ are the Hamiltonians corresponding to $\Phi$ and $\Phi^\nu$ and $\xi$ is an arbitrary but fixed boundary condition. Let us start by considering the infinite-volume limit of the first summand in \eqref{BeforeLimit2}.
We show that for a general translation-invariant IPS $L$ obeying  welldefinedness conditions as in \cite{Li85} and for $\L\uparrow\Z^d$ we have
\begin{equation}\label{Limiting 15.23}
\begin{split}
|\frac{1}{|\L|}\sum_{\o_\L}\nu(L1_{\o_\L})H_\L(\o_\L\xi_{\L^c})-\langle\nu,\Phi\rangle_L|\to0
\end{split}
\end{equation}
where 
\begin{equation}\label{GibbsRepr}
\begin{split}
\langle\nu,\Phi\rangle_L:=\int\nu(d\eta)\sum_{\D\ni0}\int 
c_\D(\eta,d\zeta_\D)\frac{1}{|\D|}\sum_{A\cap\D\neq\emptyset}[\Phi_A(\zeta_\D\eta_{\D^c})-\Phi_A(\eta)].
\end{split}
\end{equation}
Notice that $\langle\nu,\Phi\rangle_L$ becomes $\nu(LH_0)$ if the rates are just defined for single-site jumps.
In order to prove \eqref{Limiting 15.23} let us write
\begin{equation*}\label{RotationDerivativePart2a}
\begin{split}
&\frac{1}{|\L|}\sum_{\o_\L}\nu(L1_{\o_\L})H_\L(\o_\L\xi_{\L^c})\cr
&=\frac{1}{|\L|}\int\nu(d\eta)\sum_{\D\cap\L\neq\emptyset}\int c_\D(\eta,d\zeta_\D)\sum_{A\cap\L\neq\emptyset}[\Phi_A(\zeta_{\D\cap\L}\eta_{\L\setminus\D}\xi_{\L^c})
-\Phi_A(\eta_\L\xi_{\L^c})]\cr
&=\frac{1}{|\L|}\int\nu(d\eta)\sum_{\D\subset\L}\int c_\D(\eta,d\zeta_\D)\sum_{A\cap\L\neq\emptyset}[\Phi_A(\zeta_{\D}\eta_{\L\setminus\D}\xi_{\L^c})
-\Phi_A(\eta_\L\xi_{\L^c})]\cr
&\hspace{0.8cm}+\frac{1}{|\L|}\int\nu(d\eta)\sum_{\D\cap\L\neq\emptyset,\D\not\sb\L}\int c_\D(\eta,d\zeta_\D)\sum_{A\cap\L\neq\emptyset}[\Phi_A(\zeta_{\D\cap\L}\eta_{\L\setminus\D}\xi_{\L^c})
-\Phi_A(\eta_\L\xi_{\L^c})]\cr
&=:I+II.
\end{split}
\end{equation*}
On the other hand, by translation invariance the pairing can be written as
\begin{equation*}\label{Iother}
\begin{split}
&\langle\nu,\Phi\rangle_L
=\frac{1}{|\L|}\sum_{i\in\L}\int\nu(d\eta)\sum_{\D\ni i}\int \frac{1}{|\D|}c_\D(\eta,d\zeta_\D)\sum_{A\cap\D\neq\emptyset}[\Phi_A(\zeta_\D\eta_{\D^c})-\Phi_A(\eta)]\cr
&=\frac{1}{|\L|}\sum_{i\in\L}\int\nu(d\eta)\sum_{\D\ni i,\D\sb\L}\int\frac{1}{|\D|}c_\D(\eta,d\zeta_\D)\sum_{A\cap\D\neq\emptyset}[\Phi_A(\zeta_\D\eta_{\D^c})-\Phi_A(\eta)]\cr
&\hspace{1cm}+\frac{1}{|\L|}\sum_{i\in\L}\int\nu(d\eta)\sum_{\D\ni i,\D\not\sb\L}\int\frac{1}{|\D|}c_\D(\eta,d\zeta_\D)\sum_{A\cap\D\neq\emptyset}[\Phi_A(\zeta_\D\eta_{\D^c})-\Phi_A(\eta)]\cr
&=\frac{1}{|\L|}\int\nu(d\eta)\sum_{\D\sb\L}\int c_\D(\eta,d\zeta_\D)\sum_{A\cap\D\neq\emptyset}[\Phi_A(\zeta_\D\eta_{\D^c})-\Phi_A(\eta)]\cr
&\hspace{1cm}+\frac{1}{|\L|}\sum_{i\in\L}\int\nu(d\eta)\sum_{\D\ni i,\D\not\sb\L}\int\frac{1}{|\D|}c_\D(\eta,d\zeta_\D)\sum_{A\cap\D\neq\emptyset}[\Phi_A(\zeta_\D\eta_{\D^c})-\Phi_A(\eta)]\cr
&=:III+IV.\cr
\end{split}
\end{equation*}
Defining $c_\D:=\sup_{\eta}c_\D(\eta,\{1,\dots,q\}^\D)$, for the bulk term $I-III$ we have the following estimate
\begin{equation}\label{Estimate Bulk}
\begin{split}
&|I-III|
=|\frac{1}{|\L|}\int\nu(d\eta)\sum_{\D\sb\L}\int c_\D(\eta,d\zeta_\D)\times\cr
&\hspace{2cm}\sum_{A\cap\D\neq\emptyset,A\not\sb\L}[\Phi_A(\zeta_\D\eta_{\D^c})-\Phi_A(\eta)-\Phi_A(\zeta_\D\eta_{\L\setminus\D}\xi_{\L^c})+\Phi_A(\eta_\L\xi_{\L^c})]|\cr
&\leq \frac{4}{|\L|}\sum_{\D\sb\L}c_\D\sum_{A\cap\D\neq\emptyset,A\not\sb\L}\Vert\Phi_A\Vert \leq \frac{4}{|\L|}\sum_{i\in\L}\sum_{\D\ni i,\D\sb\L}\frac{1}{|\D|}c_\D\sum_{j\in\D}\sum_{A\ni j,A\not\sb\L}\Vert\Phi_A\Vert\cr
&\leq \frac{4}{|\L|}\sum_{i\in\L}\sum_{\D\ni i,\D\sb\L}c_\D\sup_{j\in\D}\sum_{A\ni j,A\not\sb\L}\Vert\Phi_A\Vert\cr
&\leq \frac{4}{|\L|}\sum_{i\in\L}\sum_{\D\ni i,\D\not\sb\G+i}c_\D\sup_{j\in\D}\sum_{A\ni j}\Vert\Phi_A\Vert+\frac{4}{|\L|}\sum_{i\in\L}\sum_{\D\ni i,\D\sb\G+i}c_\D\sup_{j\in\D}\sum_{A\ni j,A\not\sb\L}\Vert\Phi_A\Vert\cr
&\leq 4\Vert\Phi\Vert\sum_{\D\ni0,\D\not\sb\G}c_\D+\frac{4}{|\L|}\sum_{i\in\L, \O+i\sb\L}\sum_{\D\ni i,\D\sb\G+i}c_\D\sup_{j\in\D}\sum_{A\ni j,A\not\sb\O+i}\Vert\Phi_A\Vert \cr
&\hspace{3.4cm}+ \frac{4}{|\L|}\sum_{i\in\L, \O+i\not\sb\L}\sum_{\D\ni i,\D\sb\G+i}c_\D\sup_{j\in\D}\sum_{A\ni j}\Vert\Phi_A\Vert\cr
&\leq 4\Vert\Phi\Vert\sum_{\D\ni0,\D\not\sb\G}c_\D+4\sum_{\D\ni 0,\D\sb\G}c_\D\sup_{j\in\G}\sum_{A\ni j,A\not\sb\O}\Vert\Phi_A\Vert\cr
&\hspace{3.4cm}+ 4\Vert\Phi\Vert\sum_{\D\ni0,\D\sb\G}c_\D|\{i\in\L:\O+i\not\subset\L\}|\cr
\end{split}
\end{equation}
which is true for any finite set of sites $\G$ and $\O$.
By the summability assumption $\sum_{\D\ni 0}c_\D<\infty$ (see (3.3) in \cite{Li85}) the volume $\G$ can be picked in such a way that the first summand is arbitrarily small. Now $\O$ can be chosen such that the second summand becomes also small. By letting $\L\uparrow\Z^d$, the third summand of \eqref{Estimate Bulk} goes to zero.

Finally we need to show, that the error terms $II$ and $IV$ also go to zero in the infinite-volume limit.
\begin{equation*}\label{ErrorIIErrorIV}
\begin{split}
II&=\frac{1}{|\L|}\sum_{i\in\L}\sum_{\D\ni i,\D\not\sb\L}\int\nu(d\eta)\int \frac{1}{|\D\cap\L|}c_\D(\eta,d\zeta_\D)\times\cr
&\hspace{3cm}\sum_{A\cap(\D\cap\L)\neq\emptyset}[\Phi_A(\eta_{\L\setminus\D}\zeta_{\D\setminus\L^c}\xi_{\L^c})-\Phi_A(\eta_\L\xi_{\L^c})]\cr
&\leq 2\frac{1}{|\L|}\sum_{i\in\L}\sum_{\D\ni i,\D\not\sb\L}c_\D\sup_{i\in\D}\sum_{A\ni i}\Vert\Phi_A\Vert= 2\sum_{A\ni 0}\Vert\Phi_A\Vert\frac{1}{|\L|}\sum_{i\in\L}\sum_{\D\ni i,\D\not\sb\L}c_\D\cr
IV&=\frac{1}{|\L|}\sum_{i\in\L}\sum_{\D\ni i,\D\not\sb\L}\int\nu(d\eta)\int\frac{1}{|\D|}c_\D(\eta,d\zeta_\D)\sum_{A\cap\D\neq\emptyset}[\Phi_A(\zeta_\D\eta_{\D^c})-\Phi_A(\eta)]\cr
&\leq 2\sum_{A\ni 0}\Vert\Phi_A\Vert\frac{1}{|\L|}\sum_{i\in\L}\sum_{\D\ni i,\D\not\sb\L}c_\D\cr
\end{split}
\end{equation*}
In both cases, again by the final part of the proof of Theorem 15.23 in \cite{Ge11}, one verifies convergence to zero for $\L\uparrow\Z^d$.

As for the second summand in \eqref{BeforeLimit2} the exact same arguments apply and hence we can write
\begin{equation}\label{EntropyLossIPS}
\begin{split}
g_L(\nu|\mu)=\langle\nu,\Phi\rangle_L-\langle\nu,\Phi^\nu\rangle_L.
\end{split}
\end{equation}
The mapping $\Phi\mapsto\langle\nu,\Phi\rangle_L$ is linear. It is also bounded since $$|\langle\nu,\Phi\rangle_L|\leq2\Vert\Phi\Vert_0\sum_{\D\ni0}c_\D\leq2\Vert\Phi\Vert\sum_{\D\ni0}c_\D$$
which is a finite number by assumption (see \cite{Li85} assumption 3.3). In particular it is Lipschitz continuous with Lipschitz constant $2\sum_{\D\ni0}c_\D$. The mapping $\nu\mapsto\langle\nu,\Phi\rangle_L$ is weakly continuous if
$$\eta\mapsto\sum_{\D\ni0}\frac{1}{|\D|}\int
c_\D(\eta,d\zeta_\D)[H_\D(\zeta_\D\eta_{\D^c})-H_\D(\eta)]=:F_{L,\Phi}(\eta)$$ is continuous. To see that this is indeed the case, notice that for all finite $\bar\D\subset\Z^d$ the map $\eta\mapsto\sum_{\D\ni0,\D\sb\bar\D}\frac{1}{|\D|}\int  c_\D(\eta,d\zeta_\D)[H_\D(\zeta_\D\eta_{\D^c})-H_\D(\eta)]$ is continuous as a finite sum of continuous function. Further this function is convergent as $\bar\D\uparrow\Z^d$ uniformly in $\eta$ since
\begin{equation*}\label{UniformlyConv}
\begin{split}
&\sup_\eta\big|\sum_{\D\ni0,\D\sb\bar\D}\frac{1}{|\D|}\int 
c_\D(\eta,d\zeta_\D)[H_\D(\zeta_\D\eta_{\D^c})-H_\D(\eta)]-F_{L,\Phi}(\eta)\big|\cr
&=\sup_\eta\big|\sum_{\D\ni0,\D\not\sb\bar\D}\frac{1}{|\D|}\int 
c_\D(\eta,d\zeta_\D)[H_\D(\zeta_\D\eta_{\D^c})-H_\D(\eta)]\big|\cr
&\leq 2\Vert\Phi\Vert\sum_{\D\ni0,\D\not\sb\bar\D}c_\D\to0.
\end{split}
\end{equation*}
In particular the mapping $(\nu,\Phi)\mapsto\langle\nu,\Phi\rangle_L$ is jointly continuous with respect to the weak topology of measures and the $\Vert\cdot\Vert$-topology on the Banach space of potentials. This finishes the proof.
$\Cox$

\bigskip
\textbf{Remark: }
Notice that in the expected value $\langle\nu,\Phi^{\nu}\rangle_L $
the behavior of the potential 
for atypical configurations w.r.t to the measure is suppressed. 
This suggests that the existence of a uniformly convergent potential could be relaxed. In this way a weakening of the notion of a Gibbsian trajectory may do the job.
 


\subsection{A representation of continuous-time entropy decay and more general continuity conditions}

There are numerous examples of IPS with trajectories that show non-Gibbsian behavior \cite{EnRu09,EnFeHoRe10,ErKu10,FeHoMa13b,FeHoMa13a,KuRe06,HoReZu13,KuOp08,KuNy07}. One very nice example is the infinite-temperature Ising dynamics investigated in \cite{EnFeHoRe02}. Here of course the $\o$-limit measure of any trajectory is the equidistribution. In this section we generalize Theorem \ref{CorAttractivity} to not exclude the possibly of
non-Gibbsian measures in trajectories of general IPS. We start with a representation of the entropy loss for IPS similar to \eqref{EntropyLossIPS}. Right away we can write
\begin{equation}\label{RelEntLoss}
\begin{split}
g_L(\nu|\mu)=g_L(\nu)+\langle\nu,\Phi\rangle_L
\end{split}
\end{equation}
where $\Phi$ is the potential for the $L$-invariant Gibbs measure $\mu$ and 
\begin{equation*}\label{Nup}
\begin{split}
g_L(\nu):=\lim_{\L\uparrow\Z^d}\frac{1}{|\L|}\sum_{\o_\L}\nu(L1_{\o_\L})\log\nu(\o_\L).
\end{split}
\end{equation*}
Let us express $g_L(\nu)$ as a single-site density similar to $-\langle\nu,\Phi^\nu\rangle_L$ but now for probability measures $\nu$ that have the much weaker property of being non-null instead ob being Gibbs measures.

\begin{defn}\label{Def_Non_Null}
We call a random field $\nu$ non-null if there exists a number $\d>0$ and a version of the single-site conditional probabilities such that $\d\leq\nu(\eta_0|\eta_{0^c})$ for $\nu$-a.a $\eta$.
\end{defn}

\noindent
\textbf{Examples: }1. Gibbs measures for absolutely summable potentials as well as almost Gibbsian measures as defined for example in \cite{EnMaSh99,MaReMo99} are non-null.

\medskip
2. Weakly Gibbsian measures in the sense of the definitions discussed in \cite{EnMaSh99,MaReMo99,KuLeRe04} (where the potentials only have to be absolutely convergent pointwise for a set of boundary conditions with full measure) are not necessarily non-null. The same holds for the class of intuitively weakly Gibbs measures as defined for example in \cite{EnVe04}.

\medskip
3. Consider the so-called \textit{weakly dependent measures} as defined in \cite{Pf02,LePfSu95}, these are slightly less general measures in the class of asymptotically decoupled measures. Weakly dependent measures have the defining property that there exists a number $\a(\L)$ such that $\lim_{\L\uparrow\Z^d}\frac{\a(\L)}{|\L|}=0$ and
\begin{equation}\label{WeakDepe}
\begin{split}
e^{-\a(\L)}\nu(A)\nu(B)\leq\nu(A\cap B)\leq e^{\a(\L)}\nu(A)\nu(B)
\end{split}
\end{equation}
for all measurable sets $A$ and $B$ where $A$ depends only on sites in $\L$ and $B$ depends only on sites in $\L^c$. If $\nu$ is a weakly dependent random field on $\Z^d$ with finite local state space which is also translation invariant, then $\nu$ is non-null. Indeed we have for the $\nu$-a.e uniquely defined regular conditional probabilities $e^{-\a(0)}\nu(\eta_0)\leq\nu(\eta_0|\eta_{0^c})$ and by the translation invariance we can define $e^{-\a(0)}\inf_{\eta_0:\nu(\eta_0)\neq0}\nu(\eta_0)=:\d>0$.

\medskip
4. Consider 
trajectories from the infinite-temperature Ising dynamics investigated in \cite{EnFeHoRe02} where  
$\nu_t({\eta_\L})=\int\nu_0(d\s)\prod_{i\in\L}\text{Pois}_{ t}(\s_i\to\eta_i)$. Clearly $\nu_t(\eta_0|\eta_{0^c})>\frac{1}{2}(1-e^{-2t})$. 

\medskip
5. Any IPS dynamics with sitewise independent jumps on a finite local state space, where the intensity matrix $M$ is irreducible is non-null.



\begin{prop}\label{EntropyLossGen}
Let $\nu$ be a translation-invariant and non-null probability measure and $L$ a well-defined translation-invariant IPS generator (in the sense of \cite{Li85}), then 
\begin{equation}\label{G_RepX}
\begin{split}
g_L(\nu)=\int\nu(d\eta)\sum_{\D\ni0}\int c(\eta,d\xi_{\D})\frac{1}{|\D|}\log\frac{\nu(\xi_\D|\eta_{\D^c})}{\nu(\eta_\D|\eta_{\D^c})}.
\end{split}
\end{equation}
\end{prop}
Notice that the r.h.s of \eqref{G_RepX} exists since
\begin{equation*}\label{Exists}
\begin{split}
\int\nu(d\eta)\sum_{\D\ni0}\int c(\eta,d\xi_{\D})\frac{1}{|\D|}|\log\frac{\nu(\xi_\D|\eta_{\D^c})}{\nu(\eta_\D|\eta_{\D^c})}|&\leq \sum_{\D\ni0} c_\D\frac{1}{|\D|}\log\frac{1}{\inf_{\nu\text{-a.a } \eta}\nu(\eta_\D|\eta_{\D^c})}\cr
&\leq\log\frac{1}{\d}\sum_{\D\ni0} c_\D<\infty\cr
\end{split}
\end{equation*}
where we used $\nu(\eta_\D|\eta_{\D^c})\geq \d^{|\D|}$ which can be verified using the chain rule for conditional measures.
Notice also, 
by the non-positivity of \eqref{RelEntLoss} the r.h.s of \eqref{G_RepX} is an element of 
$(-\infty,-\langle\nu,\Phi\rangle_L]$ where $-\langle\nu,\Phi\rangle_L\leq 2\Vert\Phi\Vert\sum_{\D\ni0}c_\D<\infty$ and $\Phi$ is the potential for the $L$-invariant Gibbs measure $\mu$.

\bigskip
\textbf{Proof of Proposition \ref{EntropyLossGen}: }
Before taking the infinite-volume limit we have 
\begin{equation*}\label{EntropyChange1}
\begin{split}
\frac{1}{|\L|}\sum_{\o_\L}\nu(L1_{\o_\L})\log\nu(\o_\L)&=\frac{1}{|\L|}\sum_{\o_\L}\int\nu(d\eta)\sum_{\D\cap\L\neq\emptyset}\int c(\eta,d\xi_\D)1_{\o_\L}(\eta)\log\frac{\nu(\o_{\L\setminus\D}\xi_{\D\cap\L})}{\nu(\o_\L)}\cr
&=\frac{1}{|\L|}\int\nu(d\eta)\sum_{\D\cap\L\neq\emptyset}\int c(\eta,d\xi_\D)\log\frac{\nu(\xi_{\D\cap\L}|\eta_{\L\setminus\D})}{\nu(\eta_{\D\cap\L}|\eta_{\L\setminus\D})}\cr
&=\frac{1}{|\L|}\sum_{i\in\L}\int\nu(d\eta)\sum_{\D\ni i}\int \frac{1}{|\D|}c(\eta,d\xi_\D)\log\frac{\nu(\xi_{\D\cap\L}|\eta_{\L\setminus\D})}{\nu(\eta_{\D\cap\L}|\eta_{\L\setminus\D})}.\cr
\end{split}
\end{equation*}
On the other hand by translation-invariance the r.h.s of \eqref{G_RepX} can be written as
\begin{equation*}\label{EntropyChange2}
\begin{split}
\frac{1}{|\L|}\sum_{i\in\L}\int\nu(d\eta)\sum_{\D\ni i}\int \frac{1}{|\D|}c(\eta,d\xi_{\D})\log\frac{\nu(\xi_\D|\eta_{\D^c})}{\nu(\eta_\D|\eta_{\D^c})}=:G_L(\nu).\cr
\end{split}
\end{equation*}
Thus the finite-volume difference can be expressed as
\begin{equation}\label{EntropyChangeAA}
\begin{split}
&G_L(\nu)-\frac{1}{|\L|}\sum_{\o_\L}\nu(L1_{\o_\L})\log\nu(1_{\o_\L})\cr
&=\frac{1}{|\L|}\sum_{i\in\L}\sum_{\D\ni i}\frac{1}{|\D|}\int\nu(d\eta)\int c(\eta,d\xi_{\D})[\log\frac{\nu(\xi_\D|\eta_{\D^c})}{\nu(\eta_\D|\eta_{\D^c})}-\log\frac{\nu(\xi_{\D\cap\L}|\eta_{\L\setminus\D})}{\nu(\eta_{\D\cap\L}|\eta_{\L\setminus\D})}].\cr
\end{split}
\end{equation}
By the martingale convergence theorem we have for all $\xi_\D$ and cofinal increasing sequences of finite volumes $\L\supset\D$
\begin{equation*}\label{Martingale}
\begin{split}
\lim_{\L\uparrow\Z^d}\nu(\xi_\D|\eta_{\L\setminus\D})=\nu(\xi_\D|\eta_{\D^c})
\end{split}
\end{equation*}
for $\nu$-a.a. $\eta$ and in $L^1$. Hence for fixed finite $\D\subset\Z^d$ and $\L\supset\D$ by the non-nullness condition
\begin{equation*}\label{Absch1}
\begin{split}
\Bigl|\int\nu(d\eta)\int c(\eta,d\xi_{\D})\log&\frac{\nu(\xi_{\D}|\eta_{\L\setminus\D})}{\nu(\xi_\D|\eta_{\D^c})}\Bigr|
\leq\int\nu(d\eta)\int c(\eta,d\xi_{\D})\frac{|\nu(\xi_{\D}|\eta_{\L\setminus\D})-\nu(\xi_\D|\eta_{\D^c})|}{\min\{\nu(\xi_\D|\eta_{\D^c}),\nu(\xi_{\D}|\eta_{\L\setminus\D})\}}\cr
&\leq\frac{1}{\d^{|\D|}}\sup_{\tilde\eta}\int c(\tilde\eta,d\xi_{\D})
\int\nu(d\eta)|\nu(\xi_{\D}|\eta_{\L\setminus\D})-\nu(\xi_\D|\eta_{\D^c})|\cr
&\leq\frac{c_\D}{\d^{|\D|}}\max_{\xi_\D}
\int\nu(d\eta)|\nu(\xi_{\D}|\eta_{\L\setminus\D})-\nu(\xi_\D|\eta_{\D^c})|\cr
\end{split}
\end{equation*}
and hence by the martingale convergence this goes to zero in the infinite-volume limit. For the second summand in \eqref{EntropyChangeAA} the same arguments apply and hence for all $\D$ we have 
\begin{equation}\label{EntropyChangeBB}
\begin{split}
\int\nu(d\eta)\int c(\eta,d\xi_{\D})[\log\frac{\nu(\xi_\D|\eta_{\D^c})}{\nu(\xi_{\D\cap\L}|\eta_{\L\setminus\D})}-\log\frac{\nu(\eta_\D|\eta_{\D^c})}{\nu(\eta_{\D\cap\L}|\eta_{\L\setminus\D})}]\xrightarrow{\L\uparrow\Z^d} 0.\cr
\end{split}
\end{equation}
For any finite volumes $\G$ and $\O$ we can split the sum in \eqref{EntropyChangeAA} and write
\begin{equation*}\label{EntropyChangeCC}
\begin{split}
&\frac{1}{|\L|}\sum_{i\in\L}\sum_{\D\ni i}\frac{1}{|\D|}\int\nu(d\eta)\int c(\eta,d\xi_{\D})[\log\frac{\nu(\xi_\D|\eta_{\D^c})}{\nu(\xi_{\D\cap\L}|\eta_{\L\setminus\D})}-\log\frac{\nu(\eta_\D|\eta_{\D^c})}{\nu(\eta_{\D\cap\L}|\eta_{\L\setminus\D})}]\cr
&=\frac{1}{|\L|}\sum_{i\in\L:\G+i\sb\L}\sum_{\D\ni i,\D\subset\O+i}\frac{1}{|\D|}\int\nu(d\eta)\int c(\eta,d\xi_{\D})[\log\frac{\nu(\xi_\D|\eta_{\D^c})}{\nu(\xi_{\D\cap\L}|\eta_{\L\setminus\D})}-\log\frac{\nu(\eta_\D|\eta_{\D^c})}{\nu(\eta_{\D\cap\L}|\eta_{\L\setminus\D})}]\cr
&\hspace{0.5cm}+\frac{1}{|\L|}\sum_{i\in\L:\G+i\subset\L}\sum_{\D\ni i,\D\not\subset\O+i}\frac{1}{|\D|}\int\nu(d\eta)\int c(\eta,d\xi_{\D})[\log\frac{\nu(\xi_\D|\eta_{\D^c})}{\nu(\xi_{\D\cap\L}|\eta_{\L\setminus\D})}-\log\frac{\nu(\eta_\D|\eta_{\D^c})}{\nu(\eta_{\D\cap\L}|\eta_{\L\setminus\D})}]\cr
&\hspace{0.5cm}+\frac{1}{|\L|}\sum_{i\in\L:\G+i\not\subset\L}\sum_{\D\ni i}\frac{1}{|\D|}\int\nu(d\eta)\int c(\eta,d\xi_{\D})[\log\frac{\nu(\xi_\D|\eta_{\D^c})}{\nu(\xi_{\D\cap\L}|\eta_{\L\setminus\D})}-\log\frac{\nu(\eta_\D|\eta_{\D^c})}{\nu(\eta_{\D\cap\L}|\eta_{\L\setminus\D})}]\cr
&=:I+II+III.
\end{split}
\end{equation*}
For the boundary term $III$ we have
\begin{equation*}\label{EntropyChangeDD}
\begin{split}
|III|&\leq\frac{1}{|\L|}\sum_{i\in\L:\G+i\not\subset\L}\log\frac{1}{\d}\sum_{\D\ni 0}c_\D\frac{|\D\cap\L|+|\D|}{|\D|}
\leq\frac{\#\{i\in\L:\G+i\not\subset\L\}}{|\L|}\log\frac{1}{\d^2}\sum_{\D\ni 0}c_\D\cr
\end{split}
\end{equation*}
which goes to zero for $\L\uparrow\Z^d$. 
For the error term arising from the truncation of the rates represented by $II$, pick $\O$ such that $\sum_{\D\ni 0,\D\not\subset\O}c_\D<\e$. As a consequence we have
\begin{equation*}\label{EntropyChangeEE}
\begin{split}
|II|&\leq\log\frac{1}{\d^2}\sum_{\D\ni 0,\D\not\subset\O}c_\D<\log\frac{1}{\d^2}\e
\end{split}
\end{equation*}
by the same estimate as for $III$. Finally for the bulk term $I$ we can pick $\G(\O)$ such that in the martingale convergence \eqref{EntropyChangeBB} we have 
\begin{equation*}\label{EntropyChangeBBbb}
\begin{split}
\sup_{\D\ni0,\D\subset\O}|\int\nu(d\eta)\int c(\eta,d\xi_{\D})[\log\frac{\nu(\xi_\D|\eta_{\D^c})}{\nu(\xi_{\D\cap\L}|\eta_{\L\setminus\D})}-\log\frac{\nu(\eta_\D|\eta_{\D^c})}{\nu(\eta_{\D\cap\L}|\eta_{\L\setminus\D})}]|<\e\cr
\end{split}
\end{equation*}
for all $\G(\O)\sb\L-i$. Hence 
\begin{equation*}\label{EntropyChangeFF}
\begin{split}
I&=\frac{1}{|\L|}\sum_{i\in\L:\G(\O)\sb\L-i}\sum_{\D\ni0,\D\subset\O}\frac{1}{|\D|}\int\nu(d\eta)\int c(\eta,d\xi_{\D})\times\cr
&\hspace{5cm}[\log\frac{\nu(\xi_\D|\eta_{\D^c})}{\nu(\xi_{\D\cap\L-i}|\eta_{\L-i\setminus\D})}-\log\frac{\nu(\eta_\D|\eta_{\D^c})}{\nu(\eta_{\D\cap\L-i}|\eta_{\L-i\setminus\D})}]\cr
&\leq\e\sum_{\D\ni0,\D\subset\O}\frac{1}{|\D|}=\text{ Const }\e.
\end{split}
\end{equation*}
This finishes the proof.
$\Cox$

\bigskip
We know that with $\lim_{k\to\infty}\nu_{t_k}=\nu_*$ weakly it follows 
$$-\langle\nu_*,\Phi\rangle_L=\lim_{k\to\infty}g_L(\nu_{t_k})$$
by the continuity of $\langle\cdot,\Phi\rangle_L$. By Condition \ref{zero entropy loss condition} if $\nu_*\not\in\GG(\g^\Phi)$ we have $$-\langle\nu_*,\Phi\rangle_L> g_L(\nu_*).$$ 
Hence, in order to have the continuity result, it would be sufficient that
\begin{equation}\label{UpperSemi}
\begin{split}
\lim_{k\to\infty}g_L(\nu_{t_k})\leq g_L(\nu_*)
\end{split}
\end{equation}
which is upper semicontinuity of $g_L(\cdot)$ along the trajectory. 
Of course this semicontinuity 
may very well hold under less restrictive assumptions 
as in Theorem \ref{CorAttractivity} where we stipulate uniform Gibbsianness of the trajectory.

\medskip
For example in a situation where the potentials $\Phi^{\nu_{t_k}}$ still exist (as elements of the Banach space with norm $\Vert\cdot\Vert$) but fail to be uniformly bounded, the semicontinuity along a weakly convergent trajectory is determined by the semicontinuity of the function $\nu \mapsto \langle\nu,\Phi^{\nu}\rangle_L.$
Uniform boundedness of the $\Phi^{\nu}$'s is just a natural way to ensure continuity, but (semi-)continuity may hold even beyond such a requirement. 

\medskip
As another example take the infinite-temperature Glauber dynamics applied to an initial 
low temperature zero magnetic field Ising state in dimensions $d\geq 2$, investigated in \cite{EnFeHoRe02}.
This model shows provably non-Gibbsianness for sufficiently large times, without recovery of Gibbsianness 
along the trajectory. Still the relative entropy of the time evolved measure relative to the independent measure 
goes to zero, and the measure converges to the independent measure by elementary computations. For illustration let us use the representation from the proposition above. We have for any time-evolved starting measure $\nu_t$
$$g_L(\nu_t)=\int\nu_t(d\eta)\log\frac{\nu_t(\eta^0_0|\eta_{0^c})}{\nu_t(\eta_0|\eta_{0^c})}$$
where $\frac{\nu_t(\eta^0_0|\eta_{0^c})}{\nu_t(\eta_0|\eta_{0^c})}$ is bounded from above by $\frac{1+e^{-2t}}{1-e^{-2t}}$ and from below by $\frac{1-e^{-2t}}{1+e^{-2t}}$ and hence 
$$\int\nu_t(d\eta)\log\frac{\nu_t(\eta^0_0|\eta_{0^c})}{\nu_t(\eta_0|\eta_{0^c})}\to 0$$
as $t\uparrow\infty$. But this is continuity of $g_L$ at the equidistribution along the trajectory since the unique limiting measure is the equidistribution where of course $g_L(eq)=0.$


%

\bigskip
In general it would be nice to find conditions on a weakly convergent sequence of measures $\lim_{n\uparrow\infty}\nu_n=\nu_*$ such that 
$g_L(\nu_n)$ is upper semicontinuous. 
In the following theorem we give conditions on the finite-volume conditional probabilities of the convergent sequence of measures such that $g_L$ is continuous. 
In particular, many cases of site-independent jump-processes satisfy these conditions. 
\begin{thm}\label{MainTheorem}
Assume Condition \ref{zero entropy loss condition} holds with Gibbs measure $\mu$ for $\g^\Phi$. Further assume that
\begin{enumerate}
\item the sequence $(\nu_{t_n})_{n\in\N}$ of translation-invariant measures, propagated by some well-defined IPS $L$, converges weakly to $\nu_*$ as $t_n\uparrow\infty$,
\item for all $n\in\N$, $\nu_{t_n}$ is non-null with uniform constant $\d>0$ and
\item the martingale convergence theorem for the single-site conditional probabilities holds uniformly in $n\in\N$, more precisely for all $\xi_0\in\{1,\dots,q\}$
%
we have
\begin{equation}\label{MartingaleSingle}
\begin{split}
&\lim_{\L\uparrow\Z^d}\limsup_{n\uparrow\infty}\int\nu_{t_n}(d\eta)|\nu_{t_n}(\x_0|\eta_{\L\setminus0})-\nu_{t_n}(\x_0|\eta_{0^c})|=0.
\end{split}
\end{equation}
\end{enumerate}
Then $\nu_*$ is Gibbs for the same potential as $\mu$. 
\end{thm}
Notice that the convergence
\begin{equation*}\label{Martingale2}
\begin{split}
&\lim_{\L\uparrow\Z^d}\int\nu_{t_n}(d\eta)|\nu_{t_n}(\x_0|\eta_{\L\setminus0})-\nu_{t_n}(\x_0|\eta_{0^c})|=0
\end{split}
\end{equation*}
always holds by the martingale convergence theorem since conditional probabilities are uniformly integrable. Assumption three asks for the approach to zero to be uniform over the sequence of measures. Let us check some examples:

\bigskip
\textbf{Examples: }1. If $\nu_n\in\GG(\Phi^n)$ is a weakly convergent sequence of Gibbs measure for a sequence of potentials with uniform bound then 
\begin{equation*}\label{GibbsMartingale}
\begin{split}
&\int\nu_n(d\eta)|\nu_n(\x_0|\eta_{\L\setminus0})-\nu_n(\x_0|\eta_{0^c})|\cr
&=\int\nu_n(d\eta)\frac{\int\nu_n(d\s)|\nu_n(\x_0|\eta_{\L\setminus0}\s_{\L^c})-\nu_n(\x_0|\eta_{0^c})|\nu_n(\eta_{\L\setminus0}|\s_{\L^c})}{\int\nu_n(d\s)\nu_n(\eta_{\L\setminus0}|\s_{\L^c})}\cr
&\leq \text{Const}\sum_{A\ni0,A\not\sb\L}\Vert\Phi^n_A\Vert
\end{split}
\end{equation*}
where we used $|e^x-e^y|\leq|x-y|e^{\max\{|x|,|y|\}}$. By the convergence of the sequence of potentials (see step one in the proof of Theorem \ref{CorAttractivity}) and the uniform bound there exists a finite volume $\L$ such that $\sup_{n\in\N}\sum_{A\ni0,A\not\sb\L}\Vert\Phi^n_A\Vert<\e$. The non-nullness requirement is also satisfied by the uniform boundedness of the potentials.
%

\medskip
2. For the infinite-temperature Glauber dynamics from \cite{EnFeHoRe02} as mentioned above we have the non-nullness bound $\frac{1}{2}(1-e^{-2t_n})$ which can be bounded $t_n$-independently for $t_n\geq T$. Further
\begin{equation}\label{FlipMartingale}
\begin{split}
&\int\nu_{t_n}(d\eta)|\nu_{t_n}(\x_0|\eta_{\L\setminus0})-\nu_{t_n}(\x_0|\eta_{0^c})|\leq e^{-2t_n}\cr
\end{split}
\end{equation}
and hence there exists a finite volume $\L$ such that \eqref{FlipMartingale} becomes small uniformly in $t_n$.

\medskip
3. In order to move at least one step away from independent spin-flip dynamics to spatially dependent IPS consider the high-temperature spin-flip dynamics from \cite{EnFeHoRe02} Section 6 started either in another high-temperature Gibbs measure or in the low-non-zero-temperature $d$-dimensional Ising model. From \cite{MaNe01} we learn that
\begin{equation*}\label{DependentFlip}
\begin{split}
\nu_t(\eta_0|\eta_{0^c})=\sum_{\s_0=+,-}\int\nu_0(d\s)Z^t_0(\s_{0^c},\eta_{0^c})^{-1}e^{-H_0^t(\s_0\s_{0^c},\eta_0\eta_{0^c})}
\end{split}
\end{equation*}
where $H^t$ is a time-dependent Hamiltonian for the joint two-step distribution $\nu_0(d\s)S_t(\s,d\eta)$ and $Z^t$ the corresponding normalization. This Hamiltonian has nice locality properties collected in Theorem 6.3 in \cite{EnFeHoRe02}. In particular it is bounded uniformly also in $t$ (see formula 6.7 in \cite{EnFeHoRe02}) and hence $\nu_t$ is non-null uniformly in $t$. The uniform martingale convergence can be verified using formula 6.8 in \cite{EnFeHoRe02}.

\bigskip
The uniform $L_1$-convergence of the single-site conditional probabilities \eqref{MartingaleSingle} together with the non-nullness assumption implies the same convergence to hold for all finite-volume conditional probabilities. This is the statement of the following lemma which we use in the proof of Theorem \ref{MainTheorem}.

\begin{lem}\label{Single_to_Finite}
 Assume conditions 1,2 and 3 of Theorem \ref{MainTheorem} to hold, then 
\begin{equation*}\label{A}
\begin{split}
\lim_{\L\uparrow\Z^d}\limsup_{n\uparrow\infty}
\int\nu_{t_n}(d\eta)|\nu_{t_n}( \xi_{\D}| \eta_{\L\ba \D})-\nu_{t_n}( \xi_{\D}| \eta_{\D^c})|=0
\end{split}
\end{equation*}
for any finite-volume configuration $\xi_\D$.
\end{lem}

{\bf Proof: }
We use the fact that finite-$\D$ conditional probabilities can be expressed by single-site conditional probabilities (compare Theorem 1.33 of \cite{Ge11}) which 
allows us to get uniform convergence for finite $\D$ from the single-site condition. 
More precisely, let us begin with two sites $\D=\{ 1,2 \}$. We have that the two-site conditional probabilities 
can be expressed via one-site conditional probabilities by use of the identity 
\begin{equation*}\label{A}
\begin{split}
\nu_{t_n}( \xi_{1}\xi_{2}|\eta_{\L\ba\{1,2\}})=F_{\xi_{1}\xi_{2}}\Bigr(\bigl(\nu_{t_n}(\s_{1}|\s_2\eta_{\L\ba \{1,2\}}\bigr)_{\s_1,\s_2\in\{1,\dots,q\}^2}\Bigr)
\end{split}
\end{equation*}
where $F$ is a function from the set of $\{1,\dots,q\}\times\{1,\dots,q\}$ matrices given by 
\begin{equation*}\label{A}
\begin{split}
F_{\xi_{1}\xi_{2}}\Bigr(\bigl(a_{\s_{1},\s_2}\bigr)_{\s_1,\s_2\in\{1,\dots,q\}^2}\Bigr):=\frac{a_{\xi_{1}\xi_{2}}}{\sum_{\sigma_1\in\{1,\dots,q\}}\frac{a_{\s_{1},\xi_2}}{a_{\xi_{2},\s_1}}}.
\end{split}
\end{equation*}
By the uniform non-nullness hypothesis the matrix elements are uniformly bounded against zero by $\d>0$ and thus $F$ is uniformly continuous on the set of such matrices.
Using the {\em same} function we may also write
\begin{equation*}\label{A}
\begin{split}
\nu_{t_n}( \xi_{1}\xi_{2}|\eta_{\{1,2\}^c})=F_{\xi_{1}\xi_{2}}\Bigr(\bigl(\nu_{t_n}(\s_{1}|\s_2 \eta_{\{1,2\}^c}\bigr)_{\s_1,\s_2\in\{1,\dots,q\}^2}\Bigr).
\end{split}
\end{equation*}
Hence, for any $\e>0$ there exists a $\r>0$ such that 
\begin{equation}\label{A1}
\begin{split}
\sum_{\s_1,\s_2\in\{1,\dots,q\}^2} |\nu_{t_n}(\s_{1}|\s_2 \eta_{\{1,2\}^c})-\nu_{t_n}(\s_{1}|\s_2\eta_{\L\ba\{1,2\}})|\leq\r   
\end{split}
\end{equation}
implies that $|\nu_{t_n}(\xi_{1}\xi_{2}|\eta_{\L\ba\{1,2\}})-\nu_{t_n}(\xi_{1}\xi_{2}|\eta_{\{1,2\}^c})|\leq\e$.
From this follows that the single-site condition with one fixed spin-value in the conditioning of the form
\begin{equation}\label{A2}
\begin{split}
\lim_{\L\uparrow\Z^d}\limsup_{n\uparrow\infty}\int\nu_{t_n}(d \eta)|\nu_{t_n}(\s_{1}|\s_2 \eta_{\{1,2\}^c})-\nu_{t_n}(\s_{1}|\s_2\eta_{\L\ba\{1,2\}})|=0 
\end{split}
\end{equation}
for all $\s_1,\s_2$ implies the two-site condition 
\begin{equation*}\label{A}
\begin{split}
\lim_{\L\uparrow\Z^d}\limsup_{n\uparrow\infty}\int \nu_{t_n}(d\eta)|\nu_{t_n}(\xi_{1}\xi_2|\eta_{\{1,2\}^c})-\nu_{t_n}(\x_{1}\x_2|\eta_{\L\ba\{1,2\}})|=0 
\end{split}
\end{equation*}
for all $\xi_1,\x_2$. To see this write the last integrand as a difference of the function $F$ at the corresponding arguments 
and decompose the range of integration over the $\eta$-variable into the set where the condition \eqref{A1} holds, and the complement of this set. 

Further note that the above single-site condition \eqref{A2} itself follows from our assumption \eqref{MartingaleSingle} 
\begin{equation*}\label{A3}
\begin{split}
\lim_{\L\uparrow\Z^d}\limsup_{n\uparrow\infty}\int\nu_{t_n}(d \eta)|\nu_{t_n}(\s_{1}|\eta_{\{1\}^c})-\nu_{t_n}(\s_{1}|\eta_{\L\ba\{1\}})|=0 
\end{split}
\end{equation*}
estimating the integrand in \eqref{A2} by 
\begin{equation*}\label{A4}
\begin{split}
&|\nu_{t_n}(\s_{1}|\s_2 \eta_{\{1,2\}^c})-\nu_{t_n}(\s_{1}|\s_2\eta_{\L\ba\{1,2\}})|\cr
&\leq\frac{1}{\nu_{t_n}(\s_2|\eta_{\{2\}^c})}\sum_{\eta_2=1}^q\nu_{t_n}(\eta_2|\eta_{\{2\}^c})|\nu_{t_n}(\s_{1}|\eta_2\eta_{\{1,2\}^c})-\nu_{t_n}(\s_{1}|\eta_2\eta_{\L\ba\{1,2\}})|
\end{split}
\end{equation*}
and using for the first term on the r.h.s. the uniform non-nullness bound $\d$.

The case of general $\Delta$ follows from induction using a function analogous to the above $F$ 
to relate conditional probabilities in $\Delta$ to those in $\Delta \ba \{i\}$ and the singleton $\{i\}$. 
$\Cox$

\bigskip
\textbf{Proof of Theorem \ref{MainTheorem}: }We have for the relative entropy loss $g_L(\nu_{t_n}|\mu)=g_L(\nu_{t_n})+\langle\nu_{t_n},\Phi\rangle_L$ and need to show $\limsup_{n\uparrow\infty}|g_L(\nu_{t_n}|\mu)-g_L(\nu_{*}|\mu)|=0$ since then by Condition \ref{zero entropy loss condition} $\nu_*\in\GG(\g^\Phi)$. The energy part $\langle\cdot,\Phi\rangle_L$ is continuous and poses no problems. For the entropy part we can use the uniform non-nullness and Proposition \ref{EntropyLossGen} to write
\begin{equation*}\label{G_Rep22}
\begin{split}
|g_L(\nu_{t_n})|&=|\int\nu_{t_n}(d\eta)\sum_{\D\ni0}\int c(\eta,d\xi_{\D})\frac{1}{|\D|}\log\frac{\nu_{t_n}(\xi_\D|\eta_{\D^c})}{\nu_{t_n}(\eta_\D|\eta_{\D^c})}|\leq\log\frac{1}{\d}\sum_{\D\ni0}c_\D<\infty.
\end{split}
\end{equation*}
In order to truncate the (maybe infinite) sum, pick $\G$ such that $\log\frac{1}{\d}\sum_{\D\ni0, \D\not\sb\G}c_\D<\e/2$ and $\L\supset\G$. Let us use the following short-hand notations
 \begin{equation*}\label{Notation}
\begin{split}
l_n^{\G}(\xi,\eta)&:=\log\frac{\nu_{t_n}(\xi_\G|\eta_{\G^c})}{\nu_{t_n}(\eta_\G|\eta_{\G^c})}\hspace{1cm}l_{n,\L}^{\G}(\xi,\eta):=\log\frac{\nu_{t_n}(\xi_\G|\eta_{\L\setminus\G})}{\nu_{t_n}(\eta_\G|\eta_{\L\setminus\G})}\cr
l^{\G}(\xi,\eta)&:=\log\frac{\nu_{*}(\xi_\G|\eta_{\G^c})}{\nu_{*}(\eta_\G|\eta_{\G^c})}\hspace{1.2cm}l_{\L}^{\G}(\xi,\eta):=\log\frac{\nu_{*}(\xi_\G|\eta_{\L\setminus\G})}{\nu_{*}(\eta_\G|\eta_{\L\setminus\G})}.\cr
\end{split}
\end{equation*}
We can estimate the entropy difference 
\begin{equation*}\label{G_Rep}
\begin{split}
|g_L(\nu_{t_n})-&g_L(\nu_{*})|\cr
\leq\e&+\Bigl|\int\nu_{*}(d\eta)\sum_{\D\ni0,\D\subset\G}\int c(\eta,d\xi_{\D})\frac{1}{|\D|}[l^\D(\xi,\eta)-l_{\L}^\D(\xi,\eta)+l_{\L}^\D(\xi,\eta)]\cr
&\hspace{1cm}-\int\nu_{t_n}(d\eta)\sum_{\D\ni0,\D\subset\G}\int c(\eta,d\xi_{\D})\frac{1}{|\D|}\times\cr
&\hspace{2cm}[l_n^\D(\xi,\eta)-l_{n,\L}^\D(\xi,\eta)+l_{n,\L}^\D(\xi,\eta)-l_{\L}^\D(\xi,\eta)+l_{\L}^\D(\xi,\eta)]\Bigr|\cr
\leq\e&+\Bigl|\int\nu_{*}(d\eta)\sum_{\D\ni0,\D\subset\G}\int c(\eta,d\xi_{\D})\frac{1}{|\D|}[l^\D(\xi,\eta)-l_{\L}^\D(\xi,\eta)]\Bigr|\cr
&+\sum_{\D\ni0,\D\subset\G}c_\D\frac{1}{|\D|}\Vert l_{n,\L}^\D-l_\L^\D\Vert_\infty\cr
&+\Bigl|\int[\nu_*-\nu_{t_n}](d\eta)\sum_{\D\ni0,\D\subset\G}\int c(\eta,d\xi_{\D})\frac{1}{|\D|}[l_{\L}^\D(\xi,\eta)]\Bigr|\cr
&+\Bigl|\int\nu_{t_n}(d\eta)\sum_{\D\ni0,\D\subset\G}\int c(\eta,d\xi_{\D})\frac{1}{|\D|}[l_n^\D(\xi,\eta)-l_{n,\L}^\D(\xi,\eta)]\Bigr|\cr
=:\e&+A(\L)+B(n,\L)+C(n,\L)+D(n,\L).
%
%
\end{split}
\end{equation*}
All error terms become arbitrarily small. Indeed: For fixed $\L$, $\limsup_{n\uparrow\infty}B(n,\L)=0$ by the local convergence of the sequence of measures and the finiteness of the local state space. The same holds for $C(n,\L)$ since the sum is finite and $l_\L^\D$ are local and thus continuous functions.
For $A(\L)$ we can use martingale convergence as in \eqref{EntropyChangeBB}, more precisely we can estimate
\begin{equation*}\label{A}
\begin{split}
A(\L)&\leq\int\nu_{*}(d\eta)\sum_{\D\ni0,\D\subset\G}\int c(\eta,d\xi_{\D})\frac{1}{|\D|}\Bigl|l^\D(\xi,\eta)-l_{\L}^\D(\xi,\eta)\Bigr|\cr
&\leq\int\nu_{*}(d\eta)\sum_{\D\ni0,\D\subset\G}\int c(\eta,d\xi_{\D})\frac{1}{|\D|}\Bigl|\log\frac{\nu_{*}(\xi_\D|\eta_{\D^c})}{\nu_{*}(\xi_\D|\eta_{\L\setminus\D})}\Bigr|\cr
&\hspace{1cm}+\int\nu_{*}(d\eta)\sum_{\D\ni0,\D\subset\G}\int c(\eta,d\xi_{\D})\frac{1}{|\D|}\Bigl|\log\frac{\nu_{*}(\eta_\D|\eta_{\D^c})}{\nu_{*}(\eta_\D|\eta_{\L\setminus\D})}\Bigr|\cr
&\leq\sum_{\D\ni0,\D\subset\G} \frac{c_\D}{|\D|\d^{|\D|}}\max_{\xi_\D}\int\nu_{*}(d\eta)\Bigl|\nu_{*}(\xi_\D|\eta_{\D^c})-\nu_{*}(\xi_\D|\eta_{\L\setminus\D})\Bigr|\cr
&\hspace{1cm}+\sum_{\D\ni0,\D\subset\G}\frac{c_\D}{|\D|\d^{|\D|}}\max_{\eta_\D}\int\nu_{*}(d\eta)\Bigl|\nu_{*}(\eta_\D|\eta_{\D^c})-\nu_{*}(\eta_\D|\eta_{\L\setminus\D})\Bigr|\cr
\end{split}
\end{equation*}
which goes to zero for $\L\uparrow\Z^d$. For $D(n,\L)$ we can use the same estimate as for $A(\L)$ together with Lemma \ref{Single_to_Finite} and the fact that we can pick $\L$ large such that $\limsup_{n\uparrow\infty}D(n,\L)$ becomes small. 
$\Cox$

\bigskip
\textbf{Acknowledgement: }  
This work is supported by the Sonderforschungsbereich SFB $|$ TR12-Symmetries and Universality in Mesoscopic Systems. Christof K\"ulske thanks Universit\'{e} Paris Diderot - Paris 7 for kind hospitality and Giambattista Giacomin for stimulating discussions.

\end{document}